\documentclass[12pt]{amsart}
\usepackage{epsfig,color}

\usepackage{multicol,graphicx,floatflt}

\makeatother

\theoremstyle{definition}

\def\fnum{equation}
\newtheorem{Thm}[\fnum]{Theorem}
\newtheorem{Cor}[\fnum]{Corollary}

\newtheorem{Lem}[\fnum]{Lemma}

\numberwithin{equation}{section}

\newcommand{\nn}{{\bf{n}}}

\newcommand{\dist}{{\text {dist}}}

\def\RR{{\bold R}}
\def\SS{{\bold S}}

\newcommand{\Length}{{\text {Length}}}
\newcommand{\Energy}{{\text {Energy}}}

\newcommand{\cP}{{\mathcal{P}}}

\newcommand{\eqr}[1]{(\ref{#1})}

\begin{document}

\title[Width and mean curvature flow]{Width and mean curvature flow}

\author{Tobias H. Colding}%
\address{MIT\\
77 Massachusetts Avenue, Cambridge, MA 02139-4307\\
and Courant Institute of Mathematical Sciences\\
251 Mercer Street, New York, NY 10012.}
\author{William P. Minicozzi II}%
\address{Department of Mathematics\\
Johns Hopkins University\\
3400 N. Charles St.\\
Baltimore, MD 21218}

\thanks{The   authors
were partially supported by NSF Grants DMS  0606629 and DMS
0405695}


\email{colding@math.mit.edu  and minicozz@math.jhu.edu}

\maketitle

\section{Introduction}

 Given a Riemannian metric on the $2$-sphere, sweep the
$2$-sphere out  by a continuous one-parameter family of closed
curves starting and ending at point curves. Pull the sweepout
tight by, in a continuous way, pulling each curve as tight as
possible yet preserving the sweepout.  We show the following
useful property; see Theorem \ref{t:mm} below and cf. \cite{CM1},
\cite{CM2}, proposition 3.1 of \cite{CD}, proposition 3.1 of
\cite{Pi}, and 12.5 of \cite{Al}:

\vskip2mm
\parbox{6in}{Each curve in the tightened
sweepout whose length is close to the length of the longest curve
in the sweepout must itself be close to a closed geodesic. In
particular, there are curves in the sweepout that are close to
closed geodesics.}

\vskip2mm \noindent Finding closed geodesics on the $2$-sphere by
using sweepouts goes back to Birkhoff in 1917; see \cite{B1},
\cite{B2} and  section $2$ in \cite{Cr} about  Birkhoff's ideas.
The argument works equally well on any closed manifold, but only
produces non-trivial closed geodesics when the width, which is
defined in \eqr{e:w} below, is positive.  For instance, when $M$
is topologically a $2$-sphere, the width is loosely speaking up to
a constant the square of the length of the shortest closed curve
needed to ``pull over'' $M$. Thus Birkhoff's argument gives that
$2\pi$ times the width is realized as the length squared of a
closed geodesic.

The above useful property is virtually always implicit in any
sweepout construction of critical points for variational problems
yet it is not always recorded since most authors are only
interested in the existence of one critical point.

Similar results holds for sweepouts of manifolds by $2$-spheres
instead of circles; cf. \cite{CM2}.  The ideas are essentially the
same in the two cases, though the techniques in the curve case are
purely ad hoc whereas in the $2$-sphere case  additional
techniques, developed in the 1980s, have to be used to deal with
energy concentration (i.e., ``bubbling''); cf. \cite{Jo}.

As an application of the main result, we bound from above, by a
negative constant, the rate
of change of the width for a one-parameter family of convex
hypersurfaces that flows by mean curvature.  This estimate is
sharp and leads to a sharp estimate for the extinction time; cf.
\cite{CM1}, \cite{CM2} where a similar bound for the rate of
change for the two dimensional width is shown for homotopy
$3$-spheres evolving by the Ricci flow (see also \cite{Pe}).

\section{Existence of good sweepouts by curves}
\label{s:existe}

Let $M$ be a closed Riemannian manifold.
 Fix a large positive integer $L$ and let $\Lambda$ denote
the space of piecewise linear maps from $\SS^1$ to $M$ with
exactly $L$ breaks (possibly with unnecessary breaks) such that
the length of each geodesic segment is at most $2\pi$,
parametrized by a (constant) multiple of arclength, and with
Lipschitz bound $L$. By a linear map, we mean a (constant speed)
geodesic. Let $G \subset \Lambda$ denote the set of immersed
closed geodesics in $M$ of length at most $2\pi L$. (The energy of
a curve in $\Lambda$ is equal to its length squared divided by
$2\pi$. In other words, energy and length are essentially
equivalent.)

We will use the distance and topology on $\Lambda$ given by the
$W^{1,2}$ norm (Sobolev norm) on the space of maps from $\SS^1$ to
$M$. The simplest way to define the $W^{1,2}$ norm is to
isometrically embed the compact manifold $M$ into some Euclidean
space $\RR^N$.{\footnote{Recall that the square of the $W^{1,2}$
norm of a map $f:\SS^1 \to \RR^N$ is
$
     \int_{\SS^1}  \left( |f|^2 + |f'|^2 \right)
    $.
Thus two curves that are $W^{1,2}$ close are also $C^0$ close; cf.
\eqr{e:holder}.\label{fn:one}}} It will be convenient to scale
$\RR^N$, and thus $M$, by a constant so that it satisfies the
following: (M1) $\sup_M |A| \leq 1/16$, where $|A|^2$ is the norm
squared of the second fundamental form of $M$, i.e., the sum of
the squares of the principal curvatures (see, e.g., (1.24) on page
4 of \cite{CM3}); (M2)  the injectivity radius of $M$ is at least
$8\pi$ and the curvature is at most $1/64$, so that every geodesic
ball of radius at most $4\pi$ in $M$ is strictly geodesically
convex; (M3) if $x,y \in M$ with $|x-y|\leq 1$, then
  $\dist_M (x,y) \leq 2|x-y|$.

\subsection{The width}     \label{ss:width}

Let $\Omega$ be the set of continuous maps $\sigma :
\SS^1 \times [-1,1] \to M$ so that
 for each $t$ the map $\sigma (\cdot , t )$
is in $W^{1,2}$, the map  $t \to \sigma (\cdot , t )$
is continuous
 from $[-1,1]$ to
$W^{1,2}$, and finally
$\sigma$ maps $\SS^1 \times \{ -1 \}$ and $\SS^1 \times \{ 1 \}$ to points.
Given a map
    $\hat{\sigma} \in \Omega$, the homotopy class $\Omega_{\hat{\sigma}}$
 is defined to be the set of maps $\sigma \in \Omega$ that are homotopic to
    $\hat{\sigma}$ through maps in $\Omega$.
The width $W = W (\hat{\sigma})$ associated to the homotopy class
 $\Omega_{\hat{\sigma}}$ is defined by taking $\inf$ of $\max$ of
the energy of each slice. That is,  set
\begin{equation}    \label{e:w}
    W = \inf_{  \sigma \in \Omega_{\hat{\sigma}}  } \,
       \,  \max_{ t \in [-1, 1]} \,  \Energy \, (\sigma (\cdot , t ))
          \, ,
\end{equation}
where the energy is given by $\Energy \, (\sigma (\cdot , t )) =
\int_{\SS^1} \, \left| \partial_x \sigma (x,t) \right|^2 \, dx$.
The width is always non-negative and is positive if $\hat{\sigma}$
is in a non-trivial homotopy class.{\footnote{A particularly
interesting example is when $M$ is a topological $2$-sphere and
the induced map from $\SS^2$ to $M$ has degree one. In this case,
the width is positive and realized by a non-trivial closed
geodesic.  To see that the width is positive on non-trivial
homotopy classes, observe that if the maximal energy of a slice is
sufficiently small, then each curve $\sigma (\cdot , t)$ is
contained in a convex geodesic ball in $M$. Hence,   a geodesic
homotopy connects $\sigma$ to a path of point curves, so $\sigma$
is homotopically trivial.}}

The main theorem, Theorem \ref{t:mm},  that almost maximal slices in the
tightened sweepout are almost geodesics, is proven in subsection \ref{ss:mm}.  The proof of
this theorem as well as the
construction of the sequence of tighter and tighter sweepouts
uses a curve shortening map
that is defined in the next subsection.  We also state the key properties of
the shortening map in the next subsection, but postpone their proofs
to Section \ref{s:A} and the appendices.

\begin{figure}[htbp]
    \begin{minipage}[b]{0.3\textwidth}
\includegraphics[width=2in]{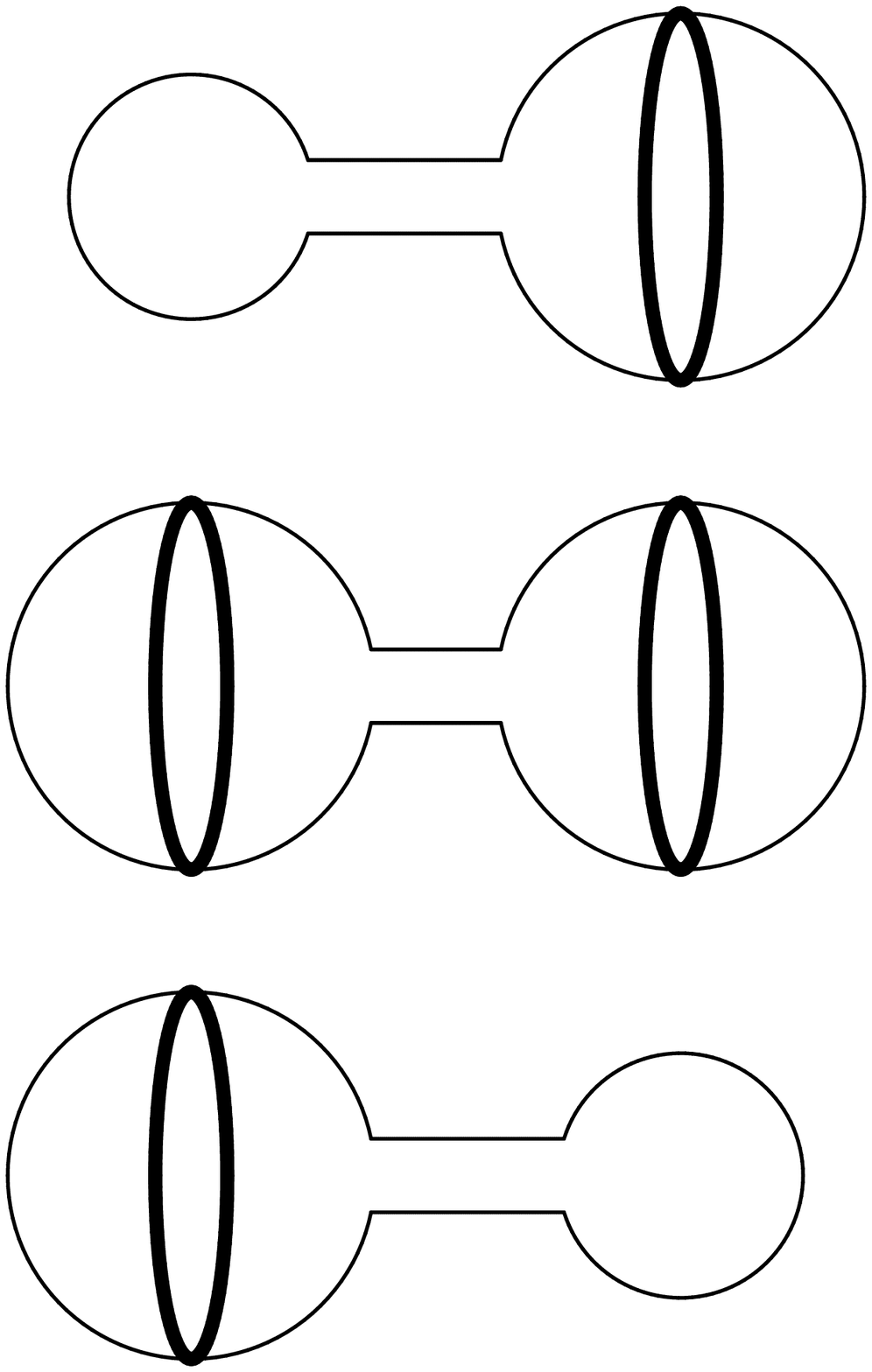}
\end{minipage}
\begin{minipage}[b]{0.5\textwidth}
\text{    }The width is continuous in the metric, but the min-max curve
that realizes it may not be.  In fact, elaborating on this example
one can easily see that the width is not in general more than
continuous in the metric.

\text{    }The continuity of the width for a
smooth one-parameter family of metrics $\{g_t\}_{t\in [0,1]}$ follows
immediately from the following:  Given $\epsilon>0$, there exists a
$\delta>0$ such that if $t\in [0,1]$ and $|s-t|<\delta$, then
$W(g_s)< W(g_t)+\epsilon$.
\end{minipage}
\end{figure}

\subsection{Curve shortening $\Psi$}    \label{ss:birk}

The curve shortening is a map $\Psi: \Lambda \to \Lambda$ so
that{\footnote{This map is essentially what is usually called
Birkhoff's curve shortening process, see section 2 of \cite{Cr}.}}
\begin{itemize}
\item[(1)] $\Psi(\gamma)$ is
homotopic to $\gamma$ and $\Length (\Psi(\gamma)) \leq \Length
(\gamma)$.
\item[(2)] $\Psi (\gamma)$ depends continuously on $\gamma$.
\item[(3)] There is a continuous function $\phi:[0,\infty) \to
[0,\infty)$ with $\phi (0) = 0$ so that
\begin{equation}    \label{e:01}
    \dist^2 (\gamma , \Psi (\gamma)) \leq \phi \left( \frac{\Length^2 (\gamma) -
\Length^2 (\Psi (\gamma))}{\Length^2 (\Psi (\gamma))} \right) \, .
\end{equation}
\item[(4)] Given $\epsilon > 0$, there exists $\delta > 0$ so that
if $\gamma \in \Lambda$ with $\dist (\gamma , G) \geq \epsilon$,
then $\Length \, (\Psi (\gamma)) \leq \Length \, (\gamma) -
\delta$.
\end{itemize}

  To define $\Psi$, we will fix a partition of   $\SS^1$
by choosing $2L$ consecutive evenly spaced points{\footnote{Note
that this is not necessarily where the piecewise linear maps have
breaks.}}
\begin{equation}
    x_0 , x_1 , x_2, \dots , x_{2L} = x_0 \in \SS^1 \, ,
\end{equation}
so that $|x_{j} - x_{j+1}| = \frac{\pi}{L}$.  $\Psi (\gamma)$ is
given in three steps.  First, we apply step 1 to $\gamma$ to get a
curve $\gamma_e$, then we apply step 2 to $\gamma_e$ to get a
curve $\gamma_o$.  In the third and final step, we reparametrize
$\gamma_o$ to get $\Psi (\gamma)$.

\vskip1mm \noindent {\bf{Step 1}}: Replace $\gamma$ on each
{\emph{even}} interval, i.e., $[x_{2j} , x_{2j+2}]$, by the linear
map with the same endpoints to get a piecewise linear curve
$\gamma_e: \SS^1 \to M$. Namely, for each $j$, we let $\gamma_e
\big|_{[x_{2j},x_{2j+2}]}$ be the unique shortest (constant speed)
geodesic from $\gamma(x_{2j})$ to $\gamma(x_{2j+2})$.

\vskip1mm \noindent {\bf{Step 2}}: Replace $\gamma_e$ on each
{\emph{odd}} interval by the linear map with the same endpoints to
get the piecewise linear curve $\gamma_o: \SS^1 \to M$.

\vskip1mm \noindent {\bf{Step 3}}: Reparametrize $\gamma_o$
(fixing $\gamma_o (x_0)$) to get the desired  constant speed
curve $\Psi(\gamma) : \SS^1 \to M$.

It is easy to see that $\Psi$ maps $\Lambda$ to $\Lambda$ and has
property (1); cf. section 2 of \cite{Cr}. Properties (2), (3) and
(4) for $\Psi$ are established in Section \ref{s:A} and Appendix
\ref{s:aB}.  Throughout the rest of this section, we will assume
these properties and use them to prove the main theorem.

\vskip2mm  The next lemma, which combines  (3) and (4), is the key
to producing the desired sequence of sweepouts.

\begin{Lem} \label{l:i1}
Given $W \geq 0$ and $\epsilon > 0$, there exists $\delta > 0$  so that
if $\gamma \in \Lambda$ and
\begin{equation}    \label{e:closei}
   2\pi \, (W - \delta) <  \Length^2 \, ( \Psi (\gamma) ) \leq  \Length^2 \, ( \gamma)  < 2\pi \, (W + \delta ) \, ,
\end{equation}
then  $\dist ( \Psi (\gamma) , G) < \epsilon$.
\end{Lem}

\begin{proof}
 If $W \leq \epsilon^2/6$,   the Wirtinger inequality{\footnote{The
 Wirtinger inequality is just the usual Poincare
inequality which bounds the $L^2$ norm in terms of the $L^2$ norm
of the derivative; i.e., $\int_{0}^{2\pi}f^2dt\leq 4\,
\int_{0}^{2\pi}(f')^2dt$ provided $f(0)=f(2\pi)=0$.
\label{fn:vert} }} gives the lemma
 with $\delta = \epsilon^2/6$.

Assume next that $W > \epsilon^2/6$.
 The triangle
inequality    gives
\begin{equation}    \label{e:closei2}
    \dist ( \Psi (\gamma) , G) \leq
        \dist ( \Psi (\gamma) , \gamma) + \dist (  \gamma , G)
     \, .
\end{equation}
Since $\Psi$  does not decrease  the length of $\gamma$ by much,
property (4) of $\Psi$ allows us to bound $\dist ( \gamma , G)$ by
$\epsilon/2$ as long as $\delta$ is sufficiently small. Similarly,
property (3) of $\Psi$ allows us to bound $\dist ( \Psi (\gamma) ,
\gamma)$ by $\epsilon/2$ as long as $\delta$ is sufficiently
small.
\end{proof}

\subsection{Defining the sweepouts}     \label{ss:sweep}

Choose a sequence of maps $\hat{\sigma}^j \in
\Omega_{\hat{\sigma}}$
 with
\begin{equation}    \label{e:minseq}
    \max_{t \in [-1,1]} \, \, \Energy \, (\hat{\sigma}^j (\cdot , t))
    < W + \frac{1}{j} \, .
\end{equation}
Observe that \eqr{e:minseq} and the Cauchy-Schwarz inequality
imply a uniform bound for the length and uniform $C^{1/2}$
continuity for the slices,   that are both independent of   $t$
and $j$. The first follows immediately and the latter follows from
\begin{align}   \label{e:holder}
      \left| \hat{\sigma}^j (x , t) \, - \,  \hat{\sigma}^j (y ,
    t)\right|^2 &\leq \left( \int_x^y \left| \partial_s \hat{\sigma}^j (s , t)
    \right| \, ds \right)^2 \notag \\
    &\leq |y-x| \,   \int_x^y \left| \partial_s \hat{\sigma}^j (s , t)
    \right|^2 \, ds  \leq |y-x| \, (W+1) \, .
\end{align}

 We will replace the $\hat{\sigma}^j$'s by sweepouts
$\sigma^j$ that, in addition to  satisfying \eqr{e:minseq}, also
satisfy that the slices $\sigma^j (\cdot , t)$ are in $\Lambda$.
  We will do this by using   local
linear replacement similar to   Step 1 of the construction of
$\Psi$.  Namely, the uniform $C^{1/2}$ bound for the slices allows
us to fix a partition of points $y_0 , \dots , y_N = y_0$ in
$\SS^1$ so that each interval $[y_i , y_{i+1}]$ is always mapped
to a ball in $M$ of radius at most  $4\pi$.  Next, for each $t$ and each $j$,
we
replace $\hat{\sigma}^j (\cdot , t) \, \big|_{[y_{i},y_{i+1}]}$ by
the linear map (geodesic) with the same endpoints  and call the
resulting map $\tilde{\sigma}^j (\cdot , t)$.  Reparametrize
$\tilde{\sigma}^j (\cdot , t)$ to have constant
speed to get $\sigma^j (\cdot , t)$.
 It is easy to see that each $\sigma^j (\cdot , t)$
satisfies \eqr{e:minseq}. Furthermore, the length bound for
$\sigma^j (\cdot , t)$ also gives a uniform Lipshitz bound for the
linear maps; let $L$ be the maximum of $N$ and this Lipshitz bound.

  It remains to show that $\sigma^j$ is continuous in the transversal
direction, i.e., with respect to $t$, and homotopic to $
\hat{\sigma}$ in $\Omega$.
 These facts were established
already by Birkhoff (see \cite{B1}, \cite{B2} and section $2$ of
\cite{Cr}), but also follow immediately from Appendix \ref{s:aB}.

Finally, applying the  replacement map $\Psi$ to   each $\sigma^j
(\cdot , t)$ gives a new  sequence of sweepouts $ \gamma^j = \Psi
(\sigma^j)$.  (By Appendix \ref{s:aB},  $\Psi$ depends
continuously on $t$ and preserves the homotopy class
$\Omega_{\hat{\sigma}}$; it is clear that $\Psi$
 fixes the constant maps at $t = \pm 1$.)

\subsection{Almost maximal implies almost critical}
\label{ss:mm}

Our main result is that this sequence $\gamma^j$ of sweepouts is tight in the
sense of the Introduction.  Namely, we have the following
theorem.

\begin{Thm}     \label{t:mm}
Given $W \geq 0$ and $\epsilon > 0$, there exist $\delta > 0$   so that
if $j
> 1/\delta$ and for some $t_0$
\begin{equation}    \label{e:close1}
    2\pi \, \Energy \, ( \gamma^j (\cdot , t_0))
        = \Length^2 \, ( \gamma^j (\cdot , t_0)) > 2\pi \, (W - \delta)  \, ,
\end{equation}
then for this $j$ we have $\dist \, \left(  \gamma^j (\cdot , t_0)
\, , \, G \right) < \epsilon$.
\end{Thm}

 \begin{proof}
 Let $\delta$ be given by
Lemma \ref{l:i1}.  By \eqr{e:close1}, \eqr{e:minseq}, and using
that $j
> 1/\delta$, we get
\begin{equation}    \label{e:close1a}
2\pi \, (W - \delta) <  \Length^2 \, ( \gamma^j (\cdot , t_0))
\leq    \Length^2 \, ( \sigma^j (\cdot , t_0)) <   2\pi \, (W+
\delta)  \, .
\end{equation}
Thus, since $\gamma^j (\cdot , t_0) = \Psi ( \sigma^j (\cdot ,
t_0))$, Lemma \ref{l:i1}   gives $\dist ( \gamma^j (\cdot , t_0)
\, , \, G) < \epsilon$, as claimed.
\end{proof}

\subsection{Parameter spaces}
Instead of using the unit interval, $[0,1]$, as the parameter
space for the circles in the sweepout and assuming that the curves
start and end in point curves, we could have used any compact set
$\cP$ and required that the curves are constant on $\partial \cP$
(or that $\partial \cP = \emptyset$).  In this case, let
$\Omega^{\cP}$ be the set of continuous maps $\sigma : \SS^1
\times \cP \to M$ so that
 for each $t \in \cP$ the curve $\sigma (\cdot , t )$
is in $W^{1,2}$, the map  $t \to \sigma (\cdot , t )$ is
continuous
 from $\cP$ to
$W^{1,2}$, and finally $\sigma$ maps $\partial \cP$  to point
curves. Given a map
    $\hat{\sigma} \in \Omega^{\cP}$, the homotopy class
$\Omega^{\cP}_{\hat{\sigma}} \subset \Omega^{\cP}$  is defined to
be the set of maps $\sigma \in \Omega^{\cP}$ that are homotopic to
    $\hat{\sigma}$ through maps in $\Omega^{\cP}$.  Finally,
    the
    width $W=W(\hat{\sigma})$ is
\begin{equation}    \label{e:wcp}
    W = \inf_{  \sigma \in \Omega^{\cP}_{\hat{\sigma}}  } \,
       \,  \max_{ t \in \cP} \,  \Energy \, (\sigma (\cdot , t ))
          \, .
\end{equation}
Theorem \ref{t:mm} holds for these general parameter spaces; the proof
is virtually the same with only
trivial changes.

\section{Rate of change of width under mean curvature flow}
\label{s:two}

Recall that a one-parameter family of smooth hypersurfaces $\{ M_t
\} \subset \RR^{n+1}$ with $n\geq 2$ {\it flows by mean curvature}
if
\begin{equation}\label{e:meanflow}
    z_t = {\bf{H}} (z) = \Delta_{M_t} z \, ,
\end{equation}
where $z$ are coordinates on $\RR^{n+1}$ and ${\bf{H}}$ is the
mean curvature vector.  By theorem $1.1$ and theorem $4.3$ in
\cite{Hu}, any smooth compact and strictly convex hypersurface in
$\RR^{n+1}$ remains smooth compact and strictly convex under the
mean curvature flow until it disappears in a point. For such a
hypersurface, the map which takes a point in $M$ to its unit
normal gives a diffeomorphism from $M$ to $\SS^n$. Since  $\SS^n =
\{ (x,y) \in \RR^2 \times \RR^{n-1} \, | \, |x|^2 + |y|^2 = 1 \}$
is equivalent to $\SS^1 \times \overline{B^{n-1}}$ where $B^{n-1}$
is the unit ball in $\RR^{n-1}$ and we collapse $\SS^1 \times \{ y
\}$ for each $y \in
\partial B^{n-1}$.
In particular,  we can fix a non-trivial homotopy class $\beta \in
\Omega^{\overline{B^{n-1}}}$ in $\pi_n(M_t)$ and define the width
$W(t) = W(\beta , M_t)$ using as parameter space
$\cP=\overline{B^{n-1}}$. It follows that the width $W(t)$ is
positive for each $t$ up until the flow $M_t$ becomes extinct.

The next is the main result of this section. It applies Theorem
\ref{t:mm} to bound the rate of change of the width $W(t)$ under
the mean curvature flow.

\begin{Thm}   \label{t:meanc}
Let $\{M_t\}_{t\geq 0}$ be a one-parameter family of smooth
compact and strictly convex hypersurfaces in $\RR^{n+1}$ flowing
by  mean curvature, then in the sense of limsup of forward
difference quotients
\begin{align} \label{e:tm1}
\frac{d}{dt}W&\leq -4\pi\, ,\\
W(t)&\leq W(0)-4\pi \, t\, .  \label{e:tm2}
\end{align}
\end{Thm}

If we have equality for $t=0$ in \eqr{e:tm1}, then for $M_0$ the
width is realized by a round circle in a plane.  Moreover, on the
circle in any direction tangent to $M_0$, but orthogonal to the
circle, the second fundamental form vanishes.  This follows from
 the cases of equality in the Cauchy-Schwarz inequality,
 the Borsuk-Fenchel inequality, and in \eqr{e:febo} below.

As a consequence of Theorem \ref{t:meanc}, we get the following
extinction result which is sharp in the case of shrinking
cylinders, where the radius of the cylinders, $r(t)$, satisfies
that $\frac{d}{dt} r^2=-2$, and, thus,
$t_{ext}=r^2(0)/2=W(0)/4\pi$.

\begin{Cor}   \label{c:meanc}
Let $\{M_t\}_{t\geq 0}$ be a one-parameter family of smooth
compact and strictly convex hypersurfaces in $\RR^{n+1}$ flowing
by mean curvature, then it becomes extinct after time at most
\begin{equation}
\frac{W(0)}{4\pi}\, .
\end{equation}
\end{Cor}

Although we have stated the results for compact convex
hypersurfaces, the arguments apply to certain types of non-compact
convex hypersurfaces; like shrinking cylinders.
The main requirement is that the ends are
``thin'' so that the width is finite.  We will not explore this
here.

\vskip1mm
 The key to proving the estimate on the rate of change of
width is the following consequence of the first variation formula
for volume (i.e., $9.3$ and $7.5$' in \cite{Si1}) and its
corollary:

\begin{Lem} \label{l:upper}
Let $M_t\subset \RR^{n+1}$ be smooth convex hypersurfaces that
flow by mean curvature.  If $\Sigma\subset M_0$ is a closed
minimal submanifold  and $\Sigma_t$ is the corresponding
submanifold in $M_t$ with volume $V_t$, then
\begin{equation}
\frac{d}{dt}_{t=0} V_t=- \int_{\Sigma} \langle {\bf{H}}_{\Sigma} ,
{\bf{H}}_{M_0} \rangle \leq -\int_{\Sigma}|{\bf{H}}_{\Sigma}|^2 \,
.\label{e:febo}
\end{equation}
Here ${\bf{H}}_{\Sigma}$ is the mean curvature vector of $\Sigma$
as a submanifold of $\RR^{n+1}$, which at $p\in \Sigma$ is equal
to the trace of the second fundamental form $A_{M_0}$ restricted
to $T_p\Sigma$ since $\Sigma$ is a minimal submanifold of $M_0$.
\end{Lem}

\begin{proof}
To get the inequality in \eqr{e:febo} we used that since $\Sigma$
is a  minimal submanifold of the convex hypersurface $M_0\subset
\RR^{n+1}$, then ${\bf{H}}_{\Sigma}$ points in the same direction
as ${\bf{H}}_{M_0}$ and  $|{\bf{H}}_{\Sigma}| \leq
|{\bf{H}}_{M_0}|$.
\end{proof}

In the first part of the next corollary, we will use the first
variation formula for the energy asserting that  if $\sigma_t
:[0,2\pi] \to \RR^{n+1}$ is a one-parameter family of curves
evolving by a vector field ${\bf{V}}$, then $\frac{d}{dt} \,
\Energy (\sigma_t) = 2\, \int_{0}^{2\pi} \langle
 \sigma_t'  ,
\nabla_{ \sigma_t' } {\bf{V}} \rangle$.

\begin{Cor} \label{c:upper}
Let $M_t$, $\Sigma$, $\Sigma_t$, ${\bf{H}}_{\Sigma}$, and $V_t$ be
as in Lemma \ref{l:upper}. If $\Sigma$ is a closed
{\underline{non-constant}} geodesic parametrized on $\SS^1$, then
$V_t$ is the length of $\Sigma_t$, ${\bf{H}}_{\Sigma_t}$ its
geodesic curvature as a curve in $\RR^{n+1}$, and
\begin{equation}
\pi \, \frac{d}{dt}_{t=0} \, \Energy (\Sigma_t) =  V_0  \,
\frac{d}{dt}_{t=0} V_t \leq -  V_0
\int_{\Sigma}|{\bf{H}}_{\Sigma}|^2
                                                    \leq -
  \left(\int_{\Sigma} |{\bf{H}}_{\Sigma}| \right)^2\leq -4\pi^2\,
.\label{e:febo1}
\end{equation}
If $\Sigma$ is a closed {\underline{non-constant}} minimal
surface, then $V_t$ is the area of $\Sigma_t$ and
\begin{equation}
\frac{d}{dt}_{t=0} V_t\leq -\int_{\Sigma}|{\bf{H}}_{\Sigma}|^2\leq
-16\pi\, \, .\label{e:febo2}
\end{equation}
\end{Cor}

\begin{proof}
The first inequality in \eqr{e:febo1} follows from Lemma
\ref{l:upper}, the second from the Cauchy-Schwarz inequality, and
the last inequality follows since by Borsuk-Fenchel's theorem
every closed curve in $\RR^{n+1}$ has total curvature at least
$2\pi$; see \cite{Bo}, \cite{Fe}.

The first inequality in \eqr{e:febo2} follows from Lemma
\ref{l:upper}.  The second inequality is $(1.4)$ in \cite{Si2},
but we include the proof.  Namely,
 use $\Delta_{\Sigma} |z|^2 = 4 + 2 \langle z ,
{\bf{H}}_{\Sigma} \rangle$ and $|\nabla_{\Sigma} |z|^2|^2 = 4
(|z|^2 -|z^{\perp}|^2)$ to compute
\begin{equation}
    \Delta_{\Sigma} \log |z|^2 = 2 \frac{\langle z , {\bf{H}}_{\Sigma}
    \rangle}{|z|^2} + 4 \frac{|z^{\perp}|^2}{|z|^4} =
    \left| \frac{1}{2}\, {\bf{H}}_{\Sigma} + 2 \, \frac{z^{\perp}}{|z|^2}
    \right|^2- \frac{1}{4} \, |{\bf{H}}_{\Sigma}|^2
     \, ,
\end{equation}
where $z$ is the position vector in $\RR^{n+1}$, and $z^{\perp}$
is the projection of $z$ to the normal space of $\Sigma$ at the
point $z$.  Applying Stokes' theorem to $-\Delta_{\Sigma} \log
|z|^2 $ gives
\begin{equation}
     \lim_{r \to 0}
     \frac{\int_{\partial B_{r} \cap \Sigma}
                       |\nabla_{\Sigma} |z|^2|} {r^2}
    \leq \frac{1}{4} \, \int_{\Sigma} |{\bf{H}}_{\Sigma}|^2 \, .
\label{e:fvsimon}
\end{equation}
Here $B_r$ is the ball of radius $r$ about $0$ in $\RR^{n+1}$.
  Since $\int_{\Sigma} |{\bf{H}}_{\Sigma}|^2$ is translation invariant, we can
translate so that $0 \in \Sigma$ and, thus, $\lim_{r \to 0}
     r^{-2} \, \int_{\partial B_{r} \cap \Sigma}
|\nabla_{\Sigma} |z|^2|$ is at least
     $4\,\pi$.
\end{proof}

The last ingredient needed in the proof of Theorem \ref{t:meanc} is the
following consequence of the first variation formula for the energy:  If
${\bf{V}}$ is a $C^2$ vector field and $\sigma_t$, $\eta_t$ are in
$W^{1,2}$, then
    \begin{equation}    \label{e:firstv}
\left| \frac{d}{dt} \,  \Energy (\eta_t) - \frac{d}{dt} \, \Energy
(\sigma_t) \right|  \leq   C \, ||{\bf{V}}||_{C^2} \, ||\sigma_t -
\eta_t ||_{W^{1,2}} \, \left( 1 + \sup |\sigma_t'|^2 \right) \, .
\end{equation}

\begin{proof}
(of Theorem \ref{t:meanc}.)
  Fix a time $\tau$.  Below ${C}$
denotes a constant depending only on $M_{\tau}$ but will be
allowed to change from inequality to inequality.  Let $\gamma^j$
be the sequence of sweepouts in   $M_{\tau}$  defined in
subsection \ref{ss:sweep}.  In particular, the maximal energy of a
slice in $\gamma^j $ goes to $W(\tau)$ as $j \to \infty$, the
$\gamma^j$'s are ``tightened'' in the sense of Theorem \ref{t:mm},
and $\gamma^j_s$ has Lipschitz bound $L$ independent of $j$ and
$s$. For $t \geq \tau$, let $\sigma_s^j(t)$ be the curve in $M_t$
that corresponds to $\gamma^j_s$  and set $e_{s,j}(t) = \Energy (
\sigma_s^j(t))$. We will use $\sigma_s^j(t)$  as a comparison to
get an upper bound for the width
 at times $t > \tau$.   The key for this is the following claim:
Given $\epsilon > 0$, there exist $\delta > 0$ and $h_0 > 0$ so
that if $j > 1/\delta$ and $0 < h < h_0$, then for all $s \in \cP$
\begin{equation}        \label{e:acomp1}
    e_{s,j}( \tau +h )
    - \max_{s_0} \,  e_{s_0,j}( \tau)  \leq
     [-4\pi + {C} \, \epsilon] \, h + {C} \, h^2 \, .
\end{equation}
 To see why
\eqr{e:acomp1} implies \eqr{e:tm1},  take the limit as $j\to
\infty$ (so that $\max_{s_0} \, e_{s_0,j}( \tau) \to W(\tau)$) in
\eqr{e:acomp1} to get
\begin{equation}        \label{e:defwq}
    \frac{W (\tau + h ) - W (\tau )}{h}
\leq    -4 \pi + {C} \, \epsilon + {C} \, h
     \, .
\end{equation}
Taking $\epsilon \to 0$ in \eqr{e:defwq} gives \eqr{e:tm1}.

  It remains to prove
\eqr{e:acomp1}. First,  let $\delta > 0$, depending on $\epsilon$
(and on $\tau$), be given by Theorem \ref{t:mm}.  Since $\beta$ is
 non-trivial in $\pi_n (M_{\tau})$, $W(\tau)$ is positive and, so, we can assume
that $\epsilon^2 < W(\tau)/3$ and $\delta < W(\tau)/3$. If $j
> 1/\delta$ and $e_{s,j}( \tau) >  W(\tau) -
\delta$, then Theorem \ref{t:mm} gives a
{\underline{non-constant}} closed geodesic $\eta$ in $M_{\tau}$
with $\dist (\eta , \gamma^j_s ) < \epsilon$. As in Lemma
\ref{l:upper}, let $\eta_t$ denote the image of $\eta$ in $M_t$.
Combining \eqr{e:febo1} and \eqr{e:firstv} with ${\bf{V}} =
{\bf{H}}_{M_t}$ and using the uniform Lipschitz bound $L$ for the
sweepouts at time $\tau$ gives
\begin{equation}        \label{e:diffAn2}
    \frac{d}{dt}_{t=\tau} e_{s,j}(t)
    \leq \frac{d}{dt}_{t=\tau} \Energy (\eta_t)
    +   {C} \, \epsilon \, \| {\bf{H}}_{M_{\tau}} \|_{C^2}
\, (1 + L^2) \,
        \leq -4 \pi   +
 {C} \, \epsilon  \, .
\end{equation}
 Since   $\sigma^j_s(t)$ is the composition of
 $\gamma^j_s$ with the smooth flow and
 $\gamma^j_s$ has  Lipschitz bound $L$ independent of $j$ and $s$, it is easy to
see that $e_{s,j}(\tau + h)$ is a smooth function of $h$ with a
uniform $C^2$ bound independent of both $j$ and $s$ near $h=0$.
  In particular, \eqr{e:diffAn2} and Taylor
expansion gives $h_0 > 0$
 (independent of $j$) so that \eqr{e:acomp1} holds
for  $s$ with $e_{s,j}( \tau) >  W(\tau) - \delta$.  In
 the remaining case, we have  $e_{s,j}( \tau) \leq  W(\tau) -
\delta$ so the continuity of $W(t)$ implies that \eqr{e:acomp1}
automatically holds after possibly shrinking $h_0 > 0$.

To get \eqr{e:tm2}, observe that for any $\epsilon
> 0$  the set $\{ t \, | \, W(t) \leq
W(0) - (4\pi - \epsilon) \, t \}$ contains $0$, is closed since
$W(t)$ is continuous,  and \eqr{e:tm1} implies that it is also
open.  Therefore, $W(t) \leq W(0) - (4\pi - \epsilon) \, t$ for
all $t$ up to the extinction time; taking $\epsilon \to 0$ gives
\eqr{e:tm2}.
\end{proof}

\subsection{$2$-Width}
Instead of defining the width by using sweepouts by closed curves,
we can define the width, $W_2$, ($2$-width) by sweeping out the
manifold by $2$-spheres, the width being the min-max value of the
energies{\footnote{The energy of a map $u:\SS^2 \to \RR^{n+1}$ is
$\frac{1}{2} \, \int_{\SS^2} |\nabla u|^2$.}} or, equivalently,
the areas of the slices in the sweepout. In \cite{CM1}, \cite{CM2}
we defined the width in this way.  Using \eqr{e:febo2} in place of
\eqr{e:febo1} and arguing much like above (cf. also with
\cite{CM1}, \cite{CM2}) we get the following (and the
corresponding extinction estimate; cf. Corollary \ref{c:meanc}):

\begin{Thm}   \label{t:meanc2}
Let $\{M_t\}_{t\geq 0}$ be a one-parameter family of smooth
compact and strictly convex hypersurfaces in $\RR^{n+1}$ flowing
by  mean curvature, then in the sense of limsup of forward
difference quotients
\begin{align}
\frac{d}{dt}W_2&\leq -16\pi\, ,\\
W_2(t)&\leq W_2(0)-16\pi \, t\, .
\end{align}
\end{Thm}

\section{Evolution by powers of mean curvature}

Suppose that $k>0$ and a one-parameter family of smooth
hypersurfaces $\{ M_t \} \subset \RR^{n+1}$ with $n\geq 2$  flows
by
\begin{equation}\label{e:meanpflow}
    z_t = |{\bf{H}}(z)|^k \, \nn(z) \, = |\Delta_{M_t}(z)|^k \, \nn(z) \, ,
\end{equation}
where $z$ are coordinates on $\RR^{n+1}$, $\nn = {\bf{H}}(z)/
|{\bf{H}}(z)|$ is the unit normal, and ${\bf{H}}$ is the mean
curvature vector.

In theorem $1.1$ of \cite{Sc}, F. Schulze extended Huisken's
result to evolution by any positive power of mean curvature.
Namely,
 if $M_0$ is compact, smooth,
and strictly convex, then the flow \eqr{e:meanpflow} is smooth and
remains convex
 until it becomes extinct.

Theorem \ref{t:meanc} and its corollary have analogs for these
more general flows.  Namely, we get a differential inequality for
the width,  $\frac{1}{1+k} \, \frac{d}{dt}_{t=0} W^{k+1} \leq -
(2\pi)^{(k+1)/2}$, that implies extinction in finite time.
 The proof relies on versions of Lemma \ref{l:upper} and Corollary
 \ref{c:upper} that are stated below.  The proofs of these
  are virtually the same as those in Section \ref{s:two}
  with the obvious changes.  In
 particular, we use H\"older's inequality in Corollary
 \ref{c:upperk} instead of Cauchy-Schwarz.

\begin{Lem} \label{l:upperk}
Let $M_t\subset \RR^{n+1}$ be smooth convex hypersurfaces that
flow by \eqr{e:meanpflow}.  If $\Sigma\subset M_0$ is a closed
minimal submanifold and $\Sigma_t$ is the corresponding
submanifold in $M_t$ with volume $V_t$, then
\begin{equation}
\frac{d}{dt}_{t=0} V_t=- \int_{\Sigma} \langle {\bf{H}}_{\Sigma} ,
|{\bf{H}}_{M_0}|^k \, \nn_{M_0} \rangle \leq
-\int_{\Sigma}|{\bf{H}}_{\Sigma}|^{1+k} \, .\label{e:febok}
\end{equation}
\end{Lem}

\begin{Cor} \label{c:upperk}
Let $M_t$, $\Sigma$, $\Sigma_t$, ${\bf{H}}_{\Sigma}$, and $V_t$ be
as in Lemma \ref{l:upper}. If $\Sigma$ is a closed
{\underline{non-constant}} geodesic parametrized on $\SS^1$, then
$V_t$ is the length of $\Sigma_t$, ${\bf{H}}_{\Sigma_t}$ its
geodesic curvature as a curve in $\RR^{n+1}$, and
\begin{equation}
  \frac{1}{1+k} \, \frac{d}{dt}_{t=0} V_t^{k+1}  =  V_0^k  \,
\frac{d}{dt}_{t=0} V_t \leq -  V_0^k
\int_{\Sigma}|{\bf{H}}_{\Sigma}|^{1+k}
                                                    \leq -
  \left(\int_{\Sigma} |{\bf{H}}_{\Sigma}| \right)^{k+1} \leq - (2\pi)^{k+1} \,
.\label{e:febo1a}
\end{equation}
\end{Cor}

\section{Establishing Properties (2), (3) and (4) for $\Psi$}  \label{s:A}

To prove  (2) and (3), it is useful to observe that there is an
equivalent, but more symmetric, way to construct $\Psi (\gamma)$
using four steps:
\begin{enumerate}
\item[($A_1$)] Follow Step 1 to get $\gamma_e$. \item[($B_1$)]
Reparametrize $\gamma_e$ (fixing the image of $x_0$) to get the
constant speed curve $\tilde{\gamma}_e$.  This reparametrization
moves the points $x_j$ to new points $\tilde{x}_j$ (i.e.,
$\gamma_e (x_j) = \tilde{\gamma}_e (\tilde{x}_j)$). \item[($A_2$)]
Do linear replacement on the odd $\tilde{x}_j$ intervals to get
$\tilde{\gamma}_o$. \item[($B_2$)] Reparametrize
$\tilde{\gamma}_o$  (fixing the image of $x_0$) to get the
constant speed curve $\Psi (\gamma)$.
\end{enumerate}
The reason that this gives the same curve is that
$\tilde{\gamma}_o$ is just a reparametrization of  ${\gamma}_o$.
We will also use that each of the four steps is energy
non-increasing.  This is obvious for  the linear replacements,
since linear maps minimize energy.    It follows  from the
Cauchy-Schwarz inequality for the reparametrizations, since for a
curve $\sigma:\SS^1 \to M$ we have
\begin{equation}
    \Length^2 (\sigma) \leq 2\pi \, \Energy (\sigma) \, ,
\end{equation}
with equality if and only if $|\sigma'| = \Length (\sigma)/(2\pi)$
almost everywhere.

\vskip1mm Using the alternative way of defining $\Psi (\gamma)$ in
four steps, we see that (3) follows from the triangle inequality
once   we bound $\dist (\gamma , \gamma_e)$ and $\dist (\gamma_e ,
\tilde{\gamma}_e)$ in terms of the decrease in length  (as well as
the analogs for steps $(A_2)$ and $(B_2)$).

The bound on $\dist (\gamma , \gamma_e)$  follows directly from
the following, see Appendix \ref{s:B} for the proof:

\begin{Lem} \label{l:sc}
There exists $C$ so that if  $I$ is an interval of length at most
$2\pi/L$, $\sigma_1: I\to M$ is a Lipschitz curve with
$|\sigma_1'| \leq L$, and $\sigma_2 : I\to M$ is the minimizing
geodesic with the same endpoints, then
\begin{equation}    \label{e:a2}
    \dist^2 (\sigma_1 , \sigma_2) \leq C \, \left( \Energy (\sigma_1) -
    \Energy (\sigma_2) \right) \, .
\end{equation}
\end{Lem}

Applying Lemma \ref{l:sc} on each of the $L$ intervals  in step
$(A_1)$, we get that
\begin{equation}    \label{e:a2vv}
    \dist^2 (\gamma , \gamma_e) \leq C \, \left( \Energy (\gamma) -
    \Energy (\gamma_e) \right) \leq \frac{C}{2\pi} \, \left( \Length^2 (\gamma) -
    \Length^2 (\Psi(\gamma)) \right) \, \, .
\end{equation}
 This gives the desired bound on $\dist (\gamma , \gamma_e)$ since
 $\Length (\Psi (\gamma)) \leq 2\pi \, L$.

\vskip1mm
 In bounding $\dist (\gamma_e , \tilde{\gamma}_e)$, we will use that
  $\gamma_e$ is just the
composition $\tilde{\gamma}_e \circ P$, where $P: \SS^1 \to \SS^1$
is a monotone piecewise linear map.{\footnote{The map $P$ is
Lipschitz, but the inverse map $P^{-1}$ may not be if $\gamma_e$
is constant on an interval.}} Using that $|\tilde{\gamma}_e'| =
\Length (\tilde{\gamma}_e)/(2\pi)$ (away from the breaks) and that
the integral of $P'$ is $2\pi$, an easy calculation gives
\begin{align}   \label{e:ta}
    \int \left( P' - 1 \right)^2 &= \int ( P')^2 - 2\pi =
\int \left( \frac{|\gamma_e'|}{|\tilde{\gamma}_e' \circ P|}
\right)^2 - 2\pi    =
      \frac{4\pi^2 }{\Length^2
(\tilde{\gamma}_e)}    \, \int |\gamma_e'|^2 - 2\pi \notag
\\ &= 2\pi \,   \frac{ \Energy
    (\gamma_e) - \Energy (\tilde{\gamma}_e)}{\Energy
    (\tilde{\gamma}_e)}    \leq
2\pi \,   \frac{ \Energy
    (\gamma) - \Energy (\Psi (\gamma))}{\Energy
    (\Psi (\gamma))}
    \, .
\end{align}
Since $\gamma_e$ and $\tilde{\gamma}_e$ agree at $x_0 = x_{2L}$,
the Wirtinger inequality (footnote
 \ref{fn:vert})
  bounds $\dist^2 (\gamma_e , \tilde{\gamma}_e)$
in terms of
\begin{equation}    \label{e:tb}
    \int \, \left| (\tilde{\gamma}_e \circ P)' - \tilde{\gamma}_e' \right|^2 \leq
    2\, \int \,
\left| (\tilde{\gamma}_e' \circ P) \, P'  - \tilde{\gamma}_e'
\circ P \right|^2 +  2\, \int \, \left| \tilde{\gamma}_e' \circ P
- \tilde{\gamma}_e' \right|^2 \, .
\end{equation}
We will bound both terms on the right hand side of \eqr{e:tb} in
terms of $\int |P'-1|^2$ and then appeal to \eqr{e:ta}. To bound
the first term, use that $|\tilde{\gamma}_e'|$ is (a constant)
$\leq L$ to get
\begin{equation}
\int \, \left| (\tilde{\gamma}_e' \circ P) \, P'  -
\tilde{\gamma}_e' \circ P \right|^2 \leq   L^2 \int |P'-1|^2 \, .
\end{equation}
To bound  the second integral, we will use that when $x$ and $y$
are points in $\SS^1$ that are {\emph{not}} separated by a break
point, then $\tilde{\gamma}_e$ is a geodesic from $x$ to $y$ and,
thus, $\tilde{\gamma}_e''$ is normal to $M$ and
$|\tilde{\gamma}_e''| \leq |\tilde{\gamma}_e'|^2 \, \sup_M |A|
\leq \frac{L^2}{16}$. Therefore, integrating $\tilde{\gamma}_e''$
from $x$ to $y$ gives
\begin{equation}
\label{e:lip}
    |\tilde{\gamma}_e' (x) -
\tilde{\gamma}_e'(y)| \leq |x-y| \, \sup |\tilde{\gamma}_e''|
\leq \frac{L^2}{16} \, |x-y| \,
 .
\end{equation}
 Divide $\SS^1$ into two sets, $S_1$ and $S_2$, where $S_1$ is the set
of points within distance $(\pi \, \int |P'-1|^2 )^{1/2}$ of a
break point for $\tilde{\gamma}_e$. Since $P(x_0) = x_0$, arguing
as in \eqr{e:holder} gives $|P(x) - x| \leq (\pi \, \int |P'-1|^2
)^{1/2}$. Thus, if $x \in S_2$, then $\tilde{\gamma}_e$ is smooth
between $x$ and $P(x)$. Consequently, \eqr{e:lip} gives
\begin{equation}    \label{e:19}
    \int_{S_2} \,
\left| \tilde{\gamma}_e' \circ P  - \tilde{\gamma}_e' \right|^2
\leq  \frac{L^4}{256} \, \int_{S_2} \,   |P(s)-s|^2 \leq
\frac{L^4}{64} \, \int  \, |P'-1|^2 \, ,
\end{equation}
where the last inequality used the Wirtinger inequality. On the
other hand,
\begin{equation}    \label{e:110}
    \int_{S_1} \,
\left| \tilde{\gamma}_e' \circ P  - \tilde{\gamma}_e' \right|^2
\leq    4\, L^2 \, \Length (S_1)  \leq  8 \, L^3 \, \left(\pi \,
\int |P'-1|^2 \right)^{1/2} \, ,
\end{equation}
completing the proof of property (3).

\vskip1mm  We show (2) in Appendix \ref{s:aB}.

\vskip1mm  To prove property (4), we will argue by contradiction.
Suppose therefore that there exist $\epsilon > 0$ and a sequence
$\gamma_j \in \Lambda$ with  $\Energy (\Psi (\gamma_j)) \geq
\Energy (\gamma_j) - 1/j$ and $\dist (\gamma_j , G) \geq \epsilon
> 0$; note that the second condition implies a positive lower
bound for $\Energy (\gamma_j)$.   Observe next that the space
$\Lambda$ is compact{\footnote{Compactness of $\Lambda$ follows
 since  $\sigma \in \Lambda$ depends continuously on
the images of the $L$ break points in the compact manifold $M$.}}
and, thus, a
 subsequence of the $\gamma_j$'s must converge to some $\gamma \in
 \Lambda$.  Since property (3)
implies that $\dist (\gamma_j , \Psi (\gamma_j)) \to 0$, the $\Psi
(\gamma_j)$'s also converge to $\gamma$.  The continuity of
$\Psi$, i.e., property (2) of $\Psi$, then implies that $\Psi
(\gamma) = \gamma$.  However, this implies that $\gamma \in G$
since the only fixed points of $\Psi$ are immersed closed
geodesics. This last fact, which was used already by Birkhoff (see
section $2$ in \cite{Cr}), follows immediately from Lemma
\ref{l:sc} and \eqr{e:ta}.  However, this would contradict that
the $\gamma_j$'s remain a fixed distance from any such closed
immersed geodesic, completing the proof of (4).

\appendix

\section{Proof of Lemma \ref{l:sc}}  \label{s:B}

We will need a simple consequence of  (M1) and (M3)  in Section
\ref{s:existe}.

\begin{Lem}     \label{l:norm}
If $x,y \in M$, then $\left| (x-y)^{\perp} \right| \leq |x-y|^2$,
where $(x-y)^{\perp}$ is the normal component to $M$ at  $y$.
\end{Lem}

\begin{proof}
If $|x-y|\geq 1$, then the claim is clear. Assume therefore that
$|x-y| < 1$ and $\alpha:[0,\ell] \to M$ is a minimizing unit speed
geodesic from $y$ to $x$ with $\ell \leq 2 \,|x-y|$. Let $V$ be
the unit   normal vector $V=(x-y)^{\perp}/|(x-y)^{\perp}|$, so
$\langle \alpha'(0) , V \rangle = 0$, and observe that
\begin{align}
|(x-y)^{\perp}| &= \int_0^{\ell} \langle \alpha' (s) , V \rangle
\, ds = \int_0^{\ell} \,   \langle   \alpha'(0) + \int_0^s \,
\alpha'' (t) \, dt   \, , V \rangle  \, ds  \leq \int_0^{\ell} \,
\int_0^s \, \left| \alpha'' (t) \right| \, dt \,
ds  \notag \\
&\leq \int_0^{\ell} \, \int_0^s \, |A (\alpha (t))| \, dt  \, ds
\leq \frac{1}{2} \, \ell^2 \, \sup_M |A| \leq |x-y|^2 \, .
\end{align}
\end{proof}

\begin{proof}
(of Lemma \ref{l:sc}).
 Integrating by parts and using  that $\sigma_1$ and $\sigma_2$
are equal on $\partial I$ gives
\begin{equation}    \label{e:tma}
    \int_{I} |\sigma_1'|^2 - \int_{I} |\sigma_2'|^2 -
     \int_{I} \left| (\sigma_1 - \sigma_2)' \right|^2
    =   - 2 \, \int_{I}   \langle (\sigma_1 - \sigma_2)  ,  \sigma_2'' \rangle \equiv \kappa \, .
\end{equation}
The lemma will follow by bounding $|\kappa|$   by $\frac{1}{2} \,
        \int_{I} \left| (\sigma_1 - \sigma_2)' \right|^2$ and appealing to Wirtinger's inequality.

 Since $\sigma_2$ is a geodesic on $M$,
  $\sigma_2''$ is normal to $M$ and
 $| \sigma_2''| \leq  |\sigma_2'|^2 \, \sup_M |A| \leq \frac{|\sigma_2'|^2}{16}$.
 Thus,  Lemma \ref{l:norm} gives
\begin{equation}    \label{e:almosttan2a}
   \left| \langle (\sigma_1 - \sigma_2)  ,  \sigma_2'' \rangle \right|
   \leq |(\sigma_1 - \sigma_2)^{\perp}| \, \frac{|\sigma_2'|^2}{16}
   \leq
     |\sigma_1 - \sigma_2|^2 \, \frac{|\sigma_2'|^2}{16}
        \, .
\end{equation}
 Integrating \eqr{e:almosttan2a}, using that  $|\sigma_2'|$ is constant with $|\sigma_2'| \,
 \Length (I)  \leq 2 \pi$,  and applying
 Wirtinger's inequality gives
\begin{equation}    \label{e:gotpsia}
    \left| \kappa \right| \leq
    \frac{|\sigma_2'|^2}{8} \,
\int_{I} |\sigma_1 - \sigma_2|^2
 \leq \frac{|\sigma_2'|^2}{8} \,
     \left( \frac{ \Length (I) }{\pi} \right)^2 \, \int_{I}
                 |(\sigma_1 -\sigma_2)'|^2
                 \leq
 \frac{1}{2} \,
        \int_{I} \left| (\sigma_1 - \sigma_2)' \right|^2 \, .
\end{equation}
\end{proof}

\section{The continuity of $\Psi$}  \label{s:aB}

\begin{Lem} \label{l:cty}
Let $\gamma : \SS^1 \to M$  be a $W^{1,2}$ map  with $\Energy
(\gamma)   \leq L$.  If $\gamma_e$ and $\tilde{\gamma}_e$ are
given by applying steps $(A_1)$ and $(B_1)$ to $\gamma$, then the
map $\gamma \to \tilde{\gamma}_e$ is continuous from $W^{1,2}$ to
$\Lambda$ equipped with the $W^{1,2}$ norm.
\end{Lem}

\begin{proof}
It follows from \eqr{e:holder} and the energy bound that $\dist_M
(\gamma (x_{2j}) , \gamma ( x_{2j+2}) ) \leq 2\pi$ for each $j$
and thus we can apply step $(A_1)$.  The lemma will   follow
easily from two observations:
\begin{enumerate}
\item[(C1)]  Since $W^{1,2}$ close curves are also $C^0$ close
(cf. footnote \ref{fn:one}), it follows that the points $\gamma_e
(x_{2j}) = \gamma (x_{2j})$ are continuous with respect to the
$W^{1,2}$ norm.
 \item[(C2)] Define
$\Gamma \subset M \times M$ by $ \Gamma = \{ (x,y) \in M \times M
\, | \, \dist_M (x,y) \leq 4\pi \} \, $, and define a map
$H:\Gamma \to C^1([0,1],M)$ by letting $H(x,y): [0,1] \to M$ be
the linear map from $x$ to $y$. Then the map $H$ is continuous on
$\Gamma$. Furthermore, the map $t \to H(x,y)(t)$ has uniformly
bounded first and second derivatives $|\partial_t H(x,y)| \leq 4
\pi$ and $|\partial_t^2 H(x,y)| \leq   \pi^2$; the second
derivative bound comes from (M1).
\end{enumerate}
To prove the lemma, suppose   that
   $\gamma^1 $ and $\gamma^2 $ are non-constant curves in
   $\Lambda$
(continuity  at the constant maps is obvious).  For $i = 1,2$ and
$j=1 , \dots , L$, let $a^i_j$ be the distance in $M$ from
$\gamma^i (x_{2j})$ to $\gamma^i (x_{2j+2})$.  Let $S^i =
\frac{1}{2\pi} \, \sum_{j=1}^L a^i_j$ be the speed of
$\tilde{\gamma}^i_e$, so that $|(\tilde{\gamma}^i_e)'| = S^i$
except at the $L$ break points. By (C1),  the $a^i_j$'s  are
continuous functions of $\gamma^i$ and, thus, so are $S^1$ and
$S^2$.  Moreover, (C1) and (C2) imply that $\gamma_e^1$ and
$\gamma_e^2$ are $C^1$-close on each interval $[x_{2j},x_{2j+2}]$.
Thus, we have shown that $\gamma \to \gamma_e$ is continuous.

To show that $\gamma_e \to \tilde{\gamma}_e$ is also continuous,
we will show that the $\tilde{\gamma}^i_e$'s are close when the
$\gamma_e^i$'s are.   Since the point $x_0 = x_{2L}$ is fixed under the reparametrization,
this will follow from applying Wirtinger's inequality
to $(\tilde{\gamma}^1_e -
\tilde{\gamma}^2_e) -  ( \tilde{\gamma}^1_e -
\tilde{\gamma}^2_e)(x_0)$
once we show that $\int_{\SS^1} |(\tilde{\gamma}^1_e -
\tilde{\gamma}^2_e)'|^2$ can be made small.

The piecewise linear curve
$\tilde{\gamma}^i_e$ is linear on the intervals
\begin{equation} \label{e:where}
    I^i_j = \left[ \frac{1}{S^i} \, \sum_{\ell < j} a^i_{\ell} \, , \,
        \frac{1}{S^i} \, \sum_{\ell \leq j} a^i_{\ell} \right] \, .
   \end{equation}
   Set $I_j = I^1_j \cap I^2_j$.
   Observe first that
    since the intervals $I^i_j$ in \eqr{e:where} depend
   continuously on  $\gamma_e^i$,  the measure of the complement
   $\SS^1 \setminus \left[ \cup_{j=1}^L I_j \right]$ can be made small, so that
\begin{equation} \label{e:ij2}
    \int_{\SS^1 \setminus \left[ \cup I_j \right]} \, \, \left| (\tilde{\gamma}_e^1 -  \tilde{\gamma}_e^2)'
    \right|^2 \leq  4\,  L^2 \, \Length \, \left( \SS^1 \setminus \left[ \cup I_j
    \right] \right)
\end{equation}
can also be made small.  We will divide the $I_j$'s into two groups, depending on the size of $a^1_j$.
  Fix some $\epsilon > 0$ and suppose first that $a^1_j < \epsilon$; by continuity, we can assume that
  $a^2_j < 2\epsilon$.  For such  a $j$, we get
  \begin{equation}
    \int_{I_j}  \left| (\tilde{\gamma}_e^1 -  \tilde{\gamma}_e^2)'
    \right|^2 \leq 2 \, \int_{I_j^1}  \left| (\tilde{\gamma}_e^1)' \right|^2 + 2 \int_{I_j^2}
    \left| (  \tilde{\gamma}_e^2)' \right|^2 \leq 2 \, L \, \left( a^1_j + a^2_j \right) \leq 6 \, \epsilon \, L \, .
      \end{equation}
      Since there are at most $L$ breaks, summing over these intervals
      contributes at most $6\epsilon \, L^2$ to the energy of $(\tilde{\gamma}_e^1 -
      \tilde{\gamma}_e^2)$.

      The last case to consider is an $I_j$ with $a^1_j \geq \epsilon$; by continuity, we can assume that
  $a^2_j \geq \epsilon/2$.   In this case, $\tilde{\gamma}_e^i$
      can be written on $I_j$ as the composition
    $\gamma_e^i \circ
   P^i_j$ where $\left|
   (P^i_j)'\right| = 2 \pi \, S^i/ (L a^i_j)$.  Furthermore, $P^1_j$ and $P^2_j$ both map $I_j$ into $[x_{2j}, x_{2j+2}]$
   and
   \begin{equation} \label{e:ij}
    \int_{I_j} \left| (\tilde{\gamma}_e^1 -  \tilde{\gamma}_e^2)'
    \right|^2 =  \int_{I_j} \left| ({\gamma}_e^1 \circ P_j^1 -  {\gamma}_e^2 \circ
    P_j^2)' \right|^2 \, .
    \end{equation}
Finally, this can be made small since the speed $\left|
   (P^i_j)'\right|$ is continuous{\footnote{The speed is continuous because of the lower bound for the $a^i_j$'s.}}    in $\gamma^i$ and the ${\gamma}_e^i$'s are
   $C^2$ bounded and $C^1$ close on $[x_{2j}, x_{2j+2}]$.
   Therefore, the integral over these intervals can also be made
   small since there are at most $L$ of them.
\end{proof}

The next result shows that $\Psi$ preserves the homotopy class of
a sweepout.

\begin{Lem}     \label{l:homotopy}
Let $\gamma \in \Omega$ satisfy $ \max_t \, \, \Energy \, (\gamma
(\cdot , t) ) \leq L$.
  If $\gamma_e$ and $\tilde{\gamma}_e$ are given by applying
steps $(A_1)$ and $(A_2)$ to each $\gamma (\cdot , t)$, then
$\gamma , \, \gamma_e$ and $\tilde{\gamma}_e$ are all homotopic in $\Omega$.
\end{Lem}

\begin{proof}
Given $x,y \in M$ with $\dist_M (x,y) \leq 4\pi$, let $H(x,y):
[0,1] \to M$ be the linear map from $x$ to $y$  as in (C2). It
follows that
\begin{equation}
    F(x,t,s) = H(\gamma (x,t) , \gamma_e (x,t)) (s)
\end{equation}
is an explicit homotopy with $F(\cdot , \cdot , 0) = \gamma$ and
 $F(\cdot , \cdot , 1) = \gamma_e$.

 For each $t$ with $\Length (\gamma_e (\cdot , t)) > 0$,
 $\gamma_e$ is given by $\gamma_e (\cdot , t)
= \tilde{\gamma}_e (\cdot , t) \circ
 P_t$ where $P_t$ is a monotone reparametrization of $\SS^1$ that
 fixes $x_0 = x_{2L}$.
 Moreover, $P_t$ is continuous by \eqr{e:ta} and $P_t$ depends
 continuously on $t$ by Lemma \ref{l:cty}.  Since
 $x \to (1-s)P_t(x) + sx$ gives a homotopy from $P_t$ to the
 identity map on $\SS^1$, we conclude that
\begin{equation}
    G(x,t,s) =  \tilde{\gamma}_e  \, ( (1-s)P_t(x) + sx  , t)
\end{equation}
is an explicit homotopy with $G(\cdot , \cdot , 0) = \gamma_e$ and
 $G(\cdot , \cdot , 1) = \tilde{\gamma}_e$.  Note that $P_t$ is
 not defined when  $\Length (\gamma_e (\cdot , t))= 0$, but the
 homotopy $G$ is.
\end{proof}


\begin{thebibliography}{A}




\bibitem[Al]{Al}
F.J. Almgren, The theory of varifolds, Mimeographed notes, Princeton, 1965.

\bibitem[B1]{B1}
G.D. Birkhoff, Dynamical systems with two degrees of freedom.
{\emph{TAMS}} 18 (1917), no. 2, 199--300.

\bibitem[B2]{B2}
G.D. Birkhoff, Dynamical systems, AMS Colloq. Publ. vol 9,
Providence, RI, 1927.


\bibitem[Bo]{Bo}
K. Borsuk, Sur la courbure totale des courbes fermées. Ann. Soc.
Polon. Math.  20  (1947), 251--265 (1948).


\bibitem[CD]{CD}
T.H. Colding and C. De Lellis, The min--max construction of
minimal surfaces, Surveys in differential geometry, Vol. 8,
Lectures on Geometry and Topology held in honor of Calabi, Lawson,
Siu, and Uhlenbeck at Harvard University, May 3--5, 2002,
Sponsored by   JDG, (2003) 75--107, math.AP/0303305.


\bibitem[CM1]{CM1}
T.H. Colding and W.P. Minicozzi II,  Estimates for the extinction
time for the Ricci flow on certain $3$--manifolds and a question
of Perelman, {\emph{JAMS}},  18  (2005),  no. 3, 561--569,
math.AP/0308090.


\bibitem[CM2]{CM2}
T.H. Colding and W.P. Minicozzi II, Width and finite extinction
time of Ricci flow, preprint.

\bibitem[CM3]{CM3}
T.H. Colding and W.P. Minicozzi II, Minimal surfaces, Courant
Lecture Notes in Math., v. 4, 1999.

\bibitem[Cr]{Cr} C.B. Croke, Area and the length of the shortest
closed geodesic, {\emph{J. Diff. Geom.}},  27  (1988),  no. 1,
1--21.

\bibitem[Fe]{Fe}
W. Fenchel, \"Uber Kr\"ummung und Windung geschlossener Raumkurven.
Math. Ann. 101 (1929), no. 1, 238--252.



 \bibitem[Hu]{Hu}
G. Huisken, Flow by mean curvature of convex surfaces into
spheres. J. Differential Geom.  20  (1984),  no. 1, 237--266.

\bibitem[Jo]{Jo}
J. Jost, Two--dimensional geometric variational problems, J. Wiley
and Sons, Chichester, N.Y.   (1991).

\bibitem[Pe]{Pe}
G. Perelman, Finite extinction time for the solutions to the Ricci
flow on certain three--manifolds,  math.DG/0307245.

\bibitem[Pi]{Pi}
J.T. Pitts, Existence and regularity of minimal surfaces on Riemannian
manfolds, Princeton University Press, Princeton, NJ;
University of Tokyo Press, Tokyo 1981.


\bibitem[Sc]{Sc}
F. Schulze,  Evolution of convex hypersurfaces by powers of the
mean curvature, {\emph{Math. Z.}}  251  (2005)  721--733.



\bibitem[Si1]{Si1} L. Simon, Lectures on geometric measure theory,
Australian National University Centre for Mathematical Analysis,
Canberra, 1983

\bibitem[Si2]{Si2} L. Simon, Existence
of surfaces minimizing the Willmore functional, {\emph{CAG}} 1
(1993)  281--326.

\end{thebibliography}
\end{document}